\documentclass[12pt]{article}
\usepackage{WideMargins2}

\def\F{{{\mathbb F}}}
\def\Z{{{\mathbb Z}}}

\def\cj#1{{{\overline{#1}}}}

\let\a=\alpha
\let\b=\beta
\let\g=\gamma

\let\x=\times

\let\ee = e
\def\ie{{ \it i.e.}}

\def\jacobi#1#2{{{\left(\frac{#1}{#2}\right)}}}

\def\calO{{{\cal O}}}

\title{\Large\bfseries  A structure theorem for finite fields}
\author{Antonia W.\ Bluher \\ National Security Agency \\ tbluher@access4less.net}
\date{\today}
\begin {document}
\fancytitle

\begin{abstract}
We present a new structure theorem for finite fields of odd order
that relates multiplicative
and additive properties in an interesting way. This theorem has several applications, including an improved understanding of Dickson
and Chebyshev polynomials and some formulas with a number-theoretic flavor. 
This paper is an abridged version of two Math ArXiv articles by the author. 
\end{abstract}

\section{The structure theorem} \label{sec:intro}

This paper is an abridged ``less-is-more'' version of two articles \cite{Permutation, Wilson-like}, and it was written with the purpose
to accompany a lecture that will be given at the Fq14 conference in Vancouver in June 2019.  
We present a new structure theorem for finite fields of odd order
that relates multiplicative
and additive properties in an interesting way. This theorem has several applications, including an improved understanding of Dickson
and Chebyshev polynomials and some formulas with a number-theoretic flavor. 
In particular, the theorem is used to prove:

\begin{enumerate}
\item[(i)]  $\prod\left\{\,a \in \F_q : \text{$a$ and $4-a$ are nonsquares}\,\right\} = 2$.
\item[(ii)] $\prod\left\{\, a \in \F_q : \text{$1-a$ and $3+a$ are nonsquares}\,\right\} =  
\begin{cases} 2 & \text{if $q=\pm 1 \pmod{12}$} \\ $-1$  & \text{otherwise.} \end{cases}$
\item[(iii)] The set $S = \left\{\, a \in \F_q : \text{$2+2a$ and $2-2a$ are squares}\,\right\}$ is preserved by Chebyshev polynomials; that is, 
$s \in S$ implies $C_k(s) \in S$, where $C_k(x)$ is the polynomial such that $C_k(\cos \theta) = \cos(k \theta)$.
(E.g. $C_2(x) = 2 x^2 - 1$ and $C_3(x) = 4 x^3 - 3 x$.)
\end{enumerate}

This article contains complete proofs of these assertions.
We assume the reader has a solid understanding of finite fields. Let $q$ be an odd prime power, 
$\F_q$ the field with $q$ elements, $\cj \F_q$ its algebraic closure, and $\F_q^\x$ its nonzero elements, considered as a group under
multiplication.
If $k>0$ is relatively prime to $q$, let $\mu_k$ denote the group of $k$th roots of unity in $\cj\F_q$. Then $\mu_k$ has cardinality $k$,
and $\mu_{q-1}=\F_q^\x$.
If $a\in\F_q$ then $\jacobi aq \in\Z$
denotes the Legendre symbol, equal to 1, $-1$, or 0 according as
$a$ is a nonzero square, a nonsquare, or~0.  This is multiplicative:
$\jacobi {ab}q = \jacobi a q \jacobi bq$.  Also, when considered as an element
of $\F_q$, $\jacobi aq = a^{(q-1)/2}$.

If $0\ne v \in \cj\F_q$, define $$\calO_v = \left\{\,v,1/v,-v,-1/v\,\right\},$$ 
which we call the {\it orbit} of $v$. If $w\in \calO_v$, then $\calO_w=\calO_v$.
Thus, the decision on which element will label the orbit is arbitrary.  
Here are some easy lemmas about $\calO_v$.  At the end of this section, the lemmas will be pieced together into one theorem.

\begin{lemma} \label{lem:wellDefined} $w \in \calO_v \iff (v+1/v)^2 = (w+1/w)^2 \iff \calO_v = \calO_w$. 
\end{lemma}

\begin{proof} The reader can verify the identity $(w+1/w)^2-(v+1/v)^2 = w^{-2}(w-v)(w+v)(w-1/v)(w+1/v)$.
The left side vanishes if and only if $(v+1/v)^2=(w+1/w)^2$, and the right side vanishes if and only if $w \in \left\{\,v,-v,1/v,-1/v\,\right\}=\calO_v$.
\end{proof}

\begin{lemma} $\#\calO_v < 4$ if and only if $v^4=1$. The orbits with fewer than 4 elements are $\left\{\,1,-1\,\right\}$ and 
$\left\{\,i,-i\,\right\}$, where $i^2=-1$.
\end{lemma}

\begin{proof} If $\#\calO_v < 4$ then two of $v,-v,1/v,-1/v$ coincide. The two elements $v,-v$ are distinct since $v$ is nonzero and the characteristic is odd.
Likewise $1/v,-1/v$ are  distinct. So the only way that two can coincide is if $v=1/v$ or $v=-1/v$. In the former case, $v^2=1$ and the orbit
is $\left\{\,1,-1\,\right\}$. In the latter case, $v^2=-1$ and the orbit is $\left\{\,i,-i\,\right\}$.
\end{proof}

\begin{lemma} \label{lem:rou} Let $f(x) = (x+1/x)^2/4$.
For $v \in \cj\F_q^\x$, $f(v) \in \F_q \iff v \in 
\mu_{2(q-1)} \cup \mu_{2(q+1)}$.
\end{lemma}

\begin{proof}  
Let $\tau = f(v)$. Then
$\tau \in \F_q \iff \tau = \tau^q \iff f(v) = f(v^q) \iff
v^q \in \calO_v = \left\{\,v,-v,1/v,-1/v\,\right\}$. 
The latter condition holds if and only if
$v^{q-1}=\pm 1$ or $v^{q+1}= \pm1$, \ie, $v \in \mu_{2(q-1)} \cup \mu_{2(q+1)}$.
\end{proof}

If $k$ is even, then $\mu_k$ is closed under reciprocal and negation, so it partitions into disjoint orbits.  
The intersection of $\mu_{2(q-1)}$ and $\mu_{2(q+1)}$ is $\mu_4$, which decomposes into the two ``short obrits'' 
$ \left\{\,1,-1\,\right\}\cup \left\{\,i,-i\,\right\}$. 
All orther orbits in $\mu_{2(q-1)}\cup\mu_{2(q+1)}$ have size~4. 
The total number of orbits in $\mu_{2(q-1)}\cup\mu_{2(q+1)}$ is 
$2+(2q-2-4)/4+(2q+2-4)/4=q$.

\begin{lemma} \label{lem:bijection}
The orbits of $\mu_{2(q-1)} \cup \mu_{2(q+1)}$ are in 
bijection with $\F_q$ by the map $\calO_v \mapsto \tau=f(v)$, where
$f(x)=(x+1/x)^2/4$.
The orbit $\left\{\,1,-1\,\right\}$ corresponds to $\tau=1$ and the orbit 
$\left\{\,i,-i\,\right\}$ corresponds to $\tau = 0$.
\end{lemma}

\begin{proof} This is immediate from Lemma~\ref{lem:wellDefined}, Lemma~\ref{lem:rou}, and the observation in the previous paragraph that
$\mu_{2(q-2)}\cup\mu_{2(q+2)}$ has exactly $q$ orbits.
\end{proof}

\begin{lemma} \label{lem:sqrttau}
If $\tau \in \cj\F_q$ then  the four elements
$\sqrt{\tau\,}+\sqrt{\tau-1\,}$, $\sqrt{\tau\,}-\sqrt{\tau-1\,}$, 
$-\sqrt{\tau\,}+\sqrt{\tau-1\,}$, $-\sqrt{\tau\,}-\sqrt{\tau-1\,}$
constitute an orbit. Here any choice can be made for the square roots of 
$\tau$ and $\tau-1$. For $A\in\left\{\,1,-1\,\right\}$ we have
$$(\sqrt{\tau\,}+\sqrt{\tau-1\,})^A = \sqrt{\tau\,}+A\sqrt{\tau-1\,}.$$
\end{lemma}

\begin{proof} From the calculation $(\sqrt{\tau\,}+\sqrt{\tau-1\,})
(\sqrt{\tau\,}-\sqrt{\tau-1\,})=\tau-(\tau-1)=1$, we see that 
$\sqrt{\tau\,}-\sqrt{\tau-1\,}$
is the reciprocal of $\sqrt{\tau\,}+\sqrt{\tau-1\,}$.  The result follows.
\end{proof}

\begin{lemma} \label{lem:sqrtFrob} If $\alpha \in \F_q$ and $A = \jacobi\alpha q$ then $(\sqrt\alpha)^q = A \sqrt\alpha$. \end{lemma}
\begin{proof} 
$(\sqrt \a)^q = \sqrt\alpha\, \left((\sqrt\a)^2\right)^{(q-1)/2} 
= \sqrt\alpha\, \alpha^{(q-1)/2}
=\sqrt\a\, \jacobi\a q = A \sqrt \alpha$.
\end{proof} 

We now present our main theorem.

\begin{theorem} {\bf (Structure Theorem; \cite[Theorem 4.1]{Wilson-like}).} 
Let $f(x) = (x+1/x)^2/4$. The map $\calO_v \mapsto f(v)$ gives a bijection 
between orbits of 
$\mu_{2(q-1)}\cup \mu_{2(q+1)}$ and elements of $\F_q$. The inverse
map is $\tau \mapsto \calO_v$, where $v = \sqrt{\tau\,}+\sqrt{\tau-1\,}$ 
(and this orbit consists of the four elements 
$\pm\sqrt{\tau\,}\pm\sqrt{\tau-1\,}$). 
The two short orbits $\left\{1,-1\right\}$ and $\left\{i,-i\right\}$ correspond to $\tau = 1$ and $\tau=0$, respectively.
If $\tau \not\in\left\{\,0,1\,\right\}$ (equivalently, $v^4 \ne 1$), 
then for $A,B \in \left\{\,1,-1\,\right\}$ we have
$$\jacobi{\tau}q=A,\ \jacobi{\tau-1}q=B \iff v^{q-AB}=A.$$ 
\end{theorem}

\begin{proof} The map $\calO_v \mapsto f(v)$ determines a well-defined 
bijection from orbits of $\mu_{2(q-1)}\cup\mu_{2(q+1)}$ onto $\F_q$ 
by~Lemma~\ref{lem:bijection}.  Let $v=\sqrt{\tau\,}+\sqrt{\tau-1\,}$, so
$1/v = \sqrt{\tau\,}-\sqrt{\tau-1\,}$. Then 
$f(v)=(v+1/v)^2/4 = (2\sqrt\tau)^2/4 = \tau$.  
This shows that $\tau \mapsto \calO_v$ is the inverse bijection.
The orbit $\calO_v$ consists of the elements 
$\pm\sqrt{\tau\,} \pm \sqrt{\tau-1\,}$ by Lemma~\ref{lem:sqrttau}.
Now we prove the last assertion. Assume 
$\tau \not\in \left\{0,1\right\}$ and let $\jacobi\tau q = A$ and 
$\jacobi{\tau-1}q = B$. 
By Lemmas~\ref{lem:sqrtFrob} and~\ref{lem:sqrttau}, 
$$v^q = (\sqrt{\tau\,})^q + (\sqrt{\tau-1\,})^q 
= A \sqrt{\tau\,} + B \sqrt{\tau-1\,}
= A(\sqrt{\tau\,}+AB\sqrt{\tau-1\,}) = A  v^{AB}.$$
This shows that $v^{q-AB} = A$, as claimed.
\end{proof}

It is worth noting that if $v \in \mu_{2(q-1)}\cup \mu_{2(q+1)}$ then either $(v^{q-1})^2=1$ or $(v^{q+1})^2 = 1$, so there are $A,B \in 
\left\{\,1,-1\,\right\}$ such that $v^{q-AB}=A$.  
If $v\in \left\{1,-1\right\}$ then $v^{q-1}=v^{q+1}=1$ and 
if $v\in \left\{\,i,-i\,\right\}$ then
$v^{q-1}=\jacobi{-1}q$ and $v^{q+1} = -\jacobi{-1}q$. Thus, for $v\in\mu_4$, 
values $A,B$ are not unique.  However, for all other
$v\in \mu_{2(q-1)} \cup \mu_{2(q+1)}$, 
$A$ and $B$ are uniquely determined from~$v$.

\bigskip
\noindent{\bf Example 1:}  Illustrate the theorem for $\F_3$.  For $q=3$, $\mu_{2(q-1)}\cup\mu_{2(q+1)} = \mu_4\cup\mu_8 = \mu_8$. 
Let $\zeta$ be a primitive eighth root of unity and set $i = \zeta^2$, so $i^2=-1$. 
The orbits are $\calO_1 = \left\{\,1,-1\,\right\}$, $\calO_i = \left\{\,i,-i\,\right\}$, and $\calO_\zeta = \left\{\,\zeta,-\zeta,1/\zeta,-1/\zeta\,\right\}
=\left\{\,\zeta,\zeta^5,\zeta^7,\zeta^3\,\right\}$.  The bijection is:
$$\calO_1 \mapsto (1/4)(1+1/1)^2 = 1; 
\qquad 1 \mapsto \{\, \pm\sqrt{1} \pm \sqrt{0}\, \} = \left\{1,-1\right\}.$$
$$\calO_i \mapsto (1/4)(i+1/i)^2 = 0;\qquad 
0 \mapsto \left\{\,\pm \sqrt{0} \pm \sqrt{-1}\,\right\} = \left\{i,-i\right\}$$
$$\calO_\zeta \mapsto (\zeta+1/\zeta)^2/4 = (i+1/i+2)/4 = 2/4=2;$$ 
$$2 \mapsto  \left\{\,\pm \sqrt 2 \pm \sqrt 1 \,\right\} = \left\{\, \pm\sqrt{-1} \pm 1\,\right \} = \left\{\,\pm i \pm 1\, \right\} \stackrel ?= \calO_\zeta$$
$$(i-1)^2 = -1+1-2i=i, \text{so $i-1=\pm\zeta$ in char.~3, showing $i-1\in\calO_\zeta$.}$$
$$(i+1) = 1/(i-1),\ \text{also in $\calO_\zeta$.}$$
We illustrate the last sentence of theorem for $\tau=2$: 
$A := \jacobi \tau q = \jacobi 2 3 = -1$; 
$B := \jacobi {\tau-1}q=\jacobi 13= 1$,
therefore $v=\zeta$ should satisfy $v^{q-AB}=A$, i.e. $v^4=-1$. 

\bigskip
\noindent{\bf Exercise:} Illustrate the theorem for $\F_5$.

\bigskip

Examples 2 and 3 use the structure theorem to derive formulas for the Legendre symbols $\jacobi 2q$ and $\jacobi {-3}q$.

\bigskip\noindent
{\bf Example 2:} Since 8 divides $2(q-1)$ or $2(q+1)$ for any $q$, $\calO_\zeta$ is always an orbit of $\mu_{2(q-1)}\cup\mu_{2(q+1)}$, where $\zeta^4=-1$.
The image of $\calO_\zeta$ is $(\zeta+1/\zeta)^2/4=1/2$.  
The inverse bijection sends $1/2$ to $\left\{\,\pm \sqrt{1/2} \pm \sqrt{-1/2}\,\right\}$.
Thus, $\left\{\,\pm\zeta,\pm1/\zeta\,\right\} = \left\{\,\pm\sqrt{1/2} \pm \sqrt{-1/2}\,\right\}$. 
Let $A = \jacobi {1/2}q=\jacobi 2 q$ and $B=\jacobi{-1/2}q=\varepsilon\jacobi 2 q$,
where $\varepsilon = \jacobi{-1}q=(-1)^{(q-1)/2}$.  Note that $\varepsilon = 1 \iff (q-1)/2$ is even $\iff q \equiv 1 \pmod 4$
and $\varepsilon = -1 \iff q \equiv 3 \pmod 4$. Thus, $q\equiv \varepsilon \mod 4$.  By the theorem, we have
$\zeta^{q-AB}=A$, which can be written as $\zeta^{q-\varepsilon}=\jacobi 2 q$, or alternatively, 
\begin{equation} (-1)^{(q-\varepsilon)/4} = \jacobi 2 q. \label{eq:jacobi2} \end{equation}
This is a known result, but it is interesting to see how it relates to the theorem.

\bigskip\noindent
{\bf Example 3:} Suppose that $3\nmid q$. Then 3 divides $2(q-1)$ or $2(q+1)$. Let $\omega$ be a primitive cube root of unity, so $\omega^2+\omega+1=0$. 
Then $\calO_\omega \mapsto (\omega+1/\omega)^2/4 = (\omega^2 + 1/\omega^2+2)/4 = (\omega^2+\omega+2)/4=(-1+2)/4=1/4$.
The inverse bijection sends $\tau=1/4$ to the orbit $\left\{\,\pm\sqrt{1/4} \pm \sqrt{-3/4}\,\right\}=\left\{\,\pm(1/2) \pm (\sqrt{-3}/2)\, \right\}$.
Let $A = \jacobi{\tau}q = 1$ and $B=\jacobi{\tau-1}q = \jacobi{-3/4}q = \jacobi{-3}q$. The structure theorem guarantees that $\omega^{q-B}=1$.
Thus, $3|(q-B)=q-\jacobi {-3}q$.  In other words, $\jacobi{-3}q=1$ if and only if $q \equiv 1 \pmod 3$.
Again, this is a known result.

\bigskip\noindent
{\bf Exercise:} 
Suppose that $5\nmid q$, and let $\zeta\in\cj\F_q$ be a 
primitive fifth root of unity.  Use the structure theorem to show
that $\calO_\zeta$ corresponds to an element of $\F_q$ iff $q \equiv \pm1
\pmod 5$. Next, show that $u=4f(\zeta)-2$ satisfies $u^2+u-1=0$, and that $f(\zeta)=(2+u)/4\in \F_q \iff \jacobi{5}q=1$. Conclude
that $\jacobi{5}q = 1 \iff q \equiv \pm 1 \pmod 5$. 

\section{Dickson and Chebyshev polynomials}

The $k$th Dickson polynomial of the first kind, $D_k(x)\in\Z[x]$, is determined by the recursion
$$\text{$D_0(x)=2$, $D_1(x)=x$, and $D_{k+2}(x)=xD_{k+1}(x)-D_k(x)$ for $k\ge0$.}$$
The first few are:
$$D_0(x)=2,\ D_1(x)=x,\ D_2(x) = x^2-2,\ D_3(x) = x^3-3x,\ D_4(x) = x^4-4x^2+2,$$
$$D_5(x) = x^5-5x^3+5x,\ D_6(x) = x^6-6x^4+9x^2-2.$$
It can be shown by induction that $D_k(x)$ is a monic polynomial of degree~$k$, and 
$$D_k(u+1/u) = u^k + 1/u^k,$$
where $u$ is an indeterminate. In fact, this functional equation determines $D_k$ uniquely and could
serve as an alternate definition.  Substituting $-u$ for $u$ in the functional equation implies the well-known fact that 
$$D_k(-x)=(-1)^k D_k(x).$$
Substituing $u=1$, $i$, or $\omega$ in the functional equation (where $\omega^2+\omega+1=0$), we find
$$D_k(2)=2,\qquad D_k(0) = 
\begin{cases} 0 & \text{if $k$ is odd} \\ 
             -2 & \text{if $k\equiv 2 \pmod 4$} \\
              2 & \text{if $k\equiv 0 \pmod 4$} \end{cases}
\qquad D_k(-1) = \begin{cases} 2 & \text{if $3|k$} \\ -1 & \text{if $3\nmid k$.} \end{cases} $$

Let $m=(q-\varepsilon)/4$, where $\varepsilon=(-1)^{(q-1)/2}$. In Eq.~(\ref{eq:jacobi2}) we observed that $m \in \Z$ and $(-1)^m = \jacobi 2 q$.
We use the structure theorem to prove a new factorization for $D_m(x)$ in $\F_q[x]$.

\begin{theorem} {\bf (Factorization of Dickson polynomials;
\cite[Theorem 8.1]{Wilson-like}).} \label{thm:dicksonFactorization} In $\F_q[x]$, we have 
\begin{equation} D_m(x) = \prod \left\{\,x-b : b \in \F_q, \jacobi{2-b}q = \jacobi{2+b}q = -1 \,\right\}.\label{dicksonFac} \end{equation}
Also, 
\begin{equation} D_m(x) = \jacobi 2q \prod \left\{\,b-x : b \in \F_q, \jacobi{2-b}q = \jacobi{2+b}q = -1 \,\right\}.\label{dicksonFac2} \end{equation}
\end{theorem}

\noindent{\bf Example:} $q=23$. Then $\varepsilon = (-1)^{(23-1)/2}=-1$, $m=(q-\varepsilon)/4 = 6$,
$D_6(x) = x^6-6x^4 + 9 x^2 - 2$. It is easy to verify that $2\pm b$ are both nonsquares iff $b\in \left\{\,3,5,8,15,18,20\,\right\}$. The theorem asserts that
$$D_6(x) \equiv (x-3)(x-5)(x-8)(x-15)(x-18)(x-20) \pmod{23}.$$

\begin{proof} To prove (\ref{dicksonFac}), we'll show: ({\it i}) if $b\in\F_q$ and $2\pm b$ are nonsquares, then $D_m(b)=0$; 
and ({\it ii}) there are $m$ such $b$'s.

Write $b=v^2+1/v^2$. Then
$$D_m(b) = D_m(v^2+1/v^2)=v^{2m}+1/v^{2m} = v^{(q-\varepsilon)/2} + v^{-(q-\varepsilon)/2}.$$
Let $\calO_v \longleftrightarrow \tau$. Then $\tau = (v+1/v)^2/4 = (b+2)/4$.  Let $A = \jacobi \tau q$ and $B = \jacobi{\tau-1}q$.
Then $A=\jacobi{b+2}q=-1$ and $B = \jacobi{-1}q\jacobi{2-b}q=-\varepsilon$. 
By the structure theorem, $v^{q-AB}=A$, \ie, $v^{q-\varepsilon}=-1$.
Then $v^{(q-\varepsilon)/2}=i$ is a square root of $-1$, and so $D_m(b) = i + 1/i = 0$.

To show there are exactly $m$ such $b$'s, the above calculation shows that 
$\tau = (b+2)/4 = (v+1/v)^2/4$ such that $v^{q-\varepsilon}=-1$.
There are exactly $q-\varepsilon$ such $v$, and if $v$ is a solution then
so is every element of its orbit. Note that
$v^4\ne 1$, as otherwise $v^{q-\varepsilon}=(v^4)^m=1^m=1$.
Thus, the solutions to $v^{q-\varepsilon}=-1$ partition into exactly
$(q-\varepsilon)/4=m$ orbits, hence $m$ such $b$'s. This proves (\ref{dicksonFac}).

To prove (\ref{dicksonFac2}),  multiply each term in the product by $-1$ to obtain:
$$D_m(x)=(-1)^m \prod\left\{\,-(x-b) : b\in\F_q, \jacobi{2-b}q	= \jacobi{2+b}q = -1 \,\right\}.$$
Then use that $(-1)^m = \jacobi 2q$, by (\ref{eq:jacobi2}).

\end{proof}

Setting $x=2$ in (\ref{dicksonFac}) gives
$$ D_m(2) = \prod\left\{\,2-b : \text{$2-b$ and $2+b$ are nonsquares}\,\right\}$$
Recall $D_k(2)=2$ for all $k$, so the left side is~2.  On the right side, set $a=2-b$, so $2+b=4-a$. We obtain the surprising result:

\begin{equation} 2 = \prod\left\{\,a \in \F_q : \text{$a$ and $4-a$ are nonsquares} \,\right\}. \label{eq:wl} \end{equation}

Likewise, setting $x=2$ in (\ref{dicksonFac2}) gives
\begin{equation} \jacobi 2q 2 = \prod\left\{\,a \in \F_q : \text{$-a$ and $4+a$ are nonsquares} \,\right\}. \label{eq:wl2} \end{equation}

\noindent {\bf Examples:}  We illustrate the formula (\ref{eq:wl}) for $q=7$ and $q=9$. \\
$q=7$: The nonsquares in $\F_7$ are 3, 5, 6. Both $a$ and $4-a$ are nonsquares iff $a\in\left\{\,5,6\,\right\}$. The product of these is 
$5 \x 6 = 2 \pmod 7$.  

\noindent
$q=9$: 
$\F_9 = \left\{\, c + d i : c,d \in \left\{\,0,1,2\,\right\} \,\right\}$, 
where $i^2=-1$. The squares are $(c+di)^2= c^2-d^2-cdi$, and the nonsquares 
are $\left\{\,c + d i : cd \ne 0\,\right \}$.
Both $a$ and $4-a$ are nonsquares if and only if $a=2\pm i$, and their product is $(2+i)(2-i)=4+1=2$. 

\section{Wilson-like theorems}

Wilson's Theorem famously states that
$\prod \F_q^\x = -1$.  We give the name ``Wilson-like theorems'' for theorems in which the product over an easily-described subset of
$\F_q^\x$ gives an easily-described result. Formulas~(\ref{eq:wl}) and~(\ref{eq:wl2}) are examples of Wilson-like theorems. We can obtain many more.

\begin{theorem} \label{thm:wl} {\bf \cite[Theorem 7.1]{Wilson-like} }
 \begin{enumerate} \item $\prod\left\{\,a \in \F_q : \text{$a$ and $4-a$ are nonsquares} \,\right\} = 2$.
\item $\prod\left\{\,a \in \F_q : \text{$-a$ and $4+a$ are nonsquares} \,\right\} = \jacobi 2q 2$.
\item $\prod\left\{\,a \in \F_q^\x : \text{$1-a$ and $a+3$ are nonsquares} \,\right\} = \begin{cases} 2 & \text{if $q \equiv \pm 1 \pmod{12}$} \\
-1 & \text{otherwise.} \end{cases}$
\item Suppose that $s,t\in\F_q^\x$ and $s^2 + t^2 = 4$. Then
$$\prod\left\{\,a \in \F_q : \text{$s^2-a$ and $t^2+a$ are nonsquares} \,\right\} 
= \jacobi{2+s}q 2 = \jacobi 2 q \jacobi{2+t}q 2.$$
\end{enumerate}
\end{theorem}

\begin{proof} The first two formulas are (\ref{eq:wl}) and (\ref{eq:wl2}), and they were obtained by setting $x=2$ in 
(\ref{dicksonFac}) and (\ref{dicksonFac2}).

For the third formula, set $x=-1$ in (\ref{dicksonFac}) and use similar arguments.  
The formula for $D_k(-1)$ was given in the previous section: it is 2 or $-1$ according
as $3|k$ or $3\nmid k$.  Note that $3|m \iff 12|q-\varepsilon \iff q \equiv \pm 1 \pmod{12}$.

Now we prove the fourth formula. Let $\tau=(t+2)/4$ correspond to 
$\calO_w$, so $t+2=(w+1/w)^2$.  Then $t=w^2+1/w^2$ and $t^2-2=w^4+1/w^4$.
Since $s \ne 0$, we know $t\not\in \left\{\,2,-2\,\right\}$, so $\tau\not \in \left\{\,0,1\,\right\}$.  By the structure theorem, it follows that $w^4\ne 1$, 
and $w^{q-AB}=A$,
where $A=\jacobi{\tau}q=\jacobi{t+2}q$ and $B=\jacobi{\tau-1}q = 
\jacobi{t-2}q$.  
Note that $AB=\jacobi{(t-2)(t+2)}q=\jacobi{-s^2}q=\jacobi{-1}q=\varepsilon$.
Therefore $w^{q-\varepsilon}=\jacobi{t+2}q$.

Setting $x=t^2-2$ in (\ref{dicksonFac2}) and then substituting $b=a+t^2-2$ gives
\begin{eqnarray*} D_m(w^4+w^{-4}) &=& \jacobi 2q \prod \left\{\,b-t^2+2 : \text{$2+b$ and $2-b$ are nonsquares in $\F_q$} \,\right\} \\
&=& \jacobi 2q \prod\left\{\,a \in \F_q : \text{$a+t^2$ and $s^2-a$ are nonsquares}\,\right\}. \end{eqnarray*}
The left side is $w^{4m} + w^{-4m} = w^{q-\varepsilon} + 1/w^{q-\varepsilon} = 2 \jacobi{t+2}q$.

We have shown that 
$$\prod\left\{\,a \in \F_q : \text{$s^2-a$ and $t^2+a$ are nonsquares} \,\right\} = \jacobi 2 q \jacobi{t+2}q 2.$$
We claim that in general if $\alpha,\beta,\gamma\in\F_q$, $\a^2+\b^2=\g^2$ and 
$\a\b\ne0$, then $\jacobi{\a+\g}q=\jacobi 2q \jacobi{\b+\g}q$. 
To see this, note that
$(\a+\b+\g)^2=\a^2+\b^2+\g^2 + 2\a\b+2\a\g+2\b\g = 2(\g^2+\a\b+\a\g+\b\g) 
= 2(\a+\g)(\b+\g)$, therefore $2(\a+\g)(\b+\g)$ is a square.
It is also nonzero, because $(\a+\g)(\b+\g)$ divides 
$(\g^2-\a^2)(\g^2-\b^2)=\b^2\a^2$ and $\a\b\ne 0$ by hypothesis.
Thus, the Legendre symbol of $2(\a+\g)(\b+\g)$ is 1,
and the claim follows. 
Applying this with $\a=s$, $\b=t$, $\g=2$ gives us that $\jacobi{s+2}q = \jacobi 2 q \jacobi{t+2}q  $, which completes the proof.
\end{proof}

\noindent{\bf Example.}  Suppose that $\jacobi 2 q = 1$; equivalently $q \equiv \pm1 \pmod 8$. Let $s=t=\sqrt 2$. Then Theorem~\ref{thm:wl} says that
$$\prod\left\{\,a \in \F_q : \text{$2-a$ and $2+a$ are nonsquares}\,\right\} =  \jacobi{2+\sqrt 2} q 2.$$
For the case $q=7$, we can take $\sqrt2=3$, and $\jacobi{2+\sqrt 2} q = \jacobi 57=-1$.  So the right side is $-2$.
On the left side, $2\pm a$ are nonsquares iff $a\in \left\{\,3,4\,\right\}$ and $3\x 4= 12 = -2$ in $\F_7$.

\section{Chebyshev polynomials} \label{sec:Chebyshev}

Chebyshev polynomials of the first kind, denoted $C_k$, are determined by the property that $\cos(k\theta) = C_k(\cos\theta)$.  The first few
are given by
$$C_1(x)=x,\ C_2(x) = 2 x^2 - 1,\ C_3(x) = 4 x^3 - 3 x.$$
They are closely related to Dickson polynomials as follows:
$$2\cos(k\theta) = e^{ik\theta} + 1/e^{ik\theta} = D_k(e^{i\theta} + 1/e^{i\theta}) = D_k(2\cos\theta)$$
therefore
$$C_k(\cos\theta) = \cos(k\theta) = (1/2) D_k(2\cos\theta)$$
\ie, 
$$C_k(X) = (1/2) D_k(2X).$$

\begin{theorem}  {\bf (\cite[Theorem 9.1]{Permutation}) }
The set 
$$S = \left\{\, a \in \F_q: \text{$2a+2$ and $2a-2$ are squares}\,\right\}$$ 
is preserved by Chebyshev polynomials; that is, 
$s \in S$ implies $C_k(s) \in S$.
\end{theorem}

\begin{proof} 
Let $a \in S$ and let $\tau = (a+1)/2 \longleftrightarrow \calO_v$. Then $\tau$ and $\tau-1$ are squares.
First assume $\tau \not \in \left\{\,0,1\,\right\}$. Then by the structure theorem, $v^4\ne 1$ and $v^{q-1}=1$, i.e. $v\in \F_q^\x$.
Then $C_k(a)=C_k(2\tau-1)=(1/2)D_k(4\tau-2)=(1/2) D_k( (v+1/v)^2 - 2) = (1/2) D_k( v^2 + 1/v^2)  = (1/2) (v^{2k} + 1/v^{2k})$.
We claim this belongs to $S$.  Indeed, if $b=(1/2)(v^{2k}+1/v^{2k})$ then $2b+2= (v^k+1/v^k)^2$ and $2b-2=(v^k-1/v^k)^2$.

This proves the result when $\tau \not \in \left\{\,0,1\,\right\}$, \ie\ when $a \not \in \left\{\,1,-1\,\right\}$. If $a=1$ then $C_k(a)=(1/2)D_k(2)=1=a \in S$.
If $a=-1\in S$ then $C_k(a)=(1/2)D_k(-2)=(-1)^k$. When $k$ is odd, this equals $a$, which is assumed to be in $S$.
When $k$ is even, this equals 1, and $1\in S$ because 4 and 0 are squares.
\end{proof}

\section{Further results}

We highlight some additional applications of the structure theorem from \cite{Permutation,Wilson-like}, omitting the proofs. 

\bigskip
\noindent{\it Factorization of Dickson polynomials.} \ 
Define $E_k\in\Z[x]$ for $k\ge 0$ by
\begin{equation} \label{Edef} \text{$E_0(x)=1$, $E_1(x)=x$, and
$E_{k+2}(x) = x E_{k+1}(x) - E_k(x)$.} \end{equation}
This is called a {\it Dickson polynomial of the second kind}, and like $D_k$, it has
been widely studied. It is well known (and can easily be shown by
induction on $k$) that $E_k$ is a monic polynomial of degree $k$, and 
\begin{equation} \label{Efunctional} E_{k-1}(u+1/u) = \frac{u^k - 1/u^k}{u-1/u}
\end{equation}

We showed in this note that if $\varepsilon = (-1)^{(q-1)/2}$ and $m=(q-\varepsilon)/4$ then
$$D_{m}(x) = \prod \left\{\, x - a : \text{$2-a$ and $2+a$ are nonsquares in $\F_q$} \,\right\}.$$
It turns out that there is an analogous factorization for $E_{m-1}$ (see \cite[Theorem 8.1]{Wilson-like}):
$$E_{m-1}(x) = \prod \left\{\, x - a : \text{$2-a$ and $2+a$ are nonzero squares in $\F_q$} \,\right\}.$$

\bigskip
\noindent{\it Oddball formulas.} \ 
The following results can be found in \cite{Permutation}.

\begin{theorem}
Let $S= \left\{\, b \in \F_q : \jacobi{2-b}q=-1,\ \jacobi{2+b}q = 1 \, \right \}$.  Then $b \mapsto b^2-2$ is a permutation of $S$, and the
inverse permutation is 
$$b\mapsto \prod\left\{\,b-a : \text{$2-a$ and $2+a$ are nonsquares} \,\right\}.$$
\end{theorem}

\begin{theorem} For any $c\in \F_q$, 
\begin{eqnarray}\label{jacc}
&&\prod\left\{\,c-a : a \in\F_q,\ \jacobi{a(a+4)}q= 1 \,\right\} + \\ \nonumber
&&\prod\left\{\,c-b : b \in\F_q,\ \jacobi{b(b+4)}q= -1 \,\right\} = \jacobi{c}q.
\end{eqnarray}
\end{theorem}

\bigskip
\noindent{\it More Wilson-like theorems.} \ 
The article \cite{Wilson-like} contains closed formulas for every product of the form
$$\prod\left\{\, a \in \F_q^\x : \jacobi{k+a}q = \varepsilon_1,\ \jacobi{\ell-a}q = \varepsilon_2 \,\right\},$$
where $\varepsilon_1,\varepsilon_2 \in \left\{\,1,-1\,\right\}$ and $k,\ell\in \F_q$. Theorem~\ref{thm:wl} contains a few formulas
of this type, but there are many more.
Here are some specific examples.

\bigskip\noindent
If $j$ is a square and $4-j$ is not, then
$\prod\left\{\, a \in \F_q^\x : \text{$j-a$ and $4-j+a$ are nonsquares}\,  \right\}$ is a square root of $j$.

\bigskip\noindent
If $4-j$ is a square and $j$ is not, then
$\prod\left\{\, a \in \F_q^\x : \text{$j-a$ and $4-j+a$ are nonsquares}\,  \right\}$ is a square root of $4-j$.

\bigskip\noindent
If $j$ and $4-j$ are nonsquares, then
$$\prod \left\{\, a \in \F_q^\x : \text{$j-a$ is a nonzero square and $4-j+a$ is a nonsquare} \,\right\}$$ is a square root of $j/(4-j)$.

Other results from \cite{Wilson-like} apply to a restricted set of $q$. For example, suppose that $\jacobi 5q = 1$
(equivalently, $q\equiv \pm1 \pmod 5$ by quadratic reciprocity or by the exercise at the end of Section~\ref{sec:intro}),  and let $r = (1-\sqrt 5)/2$, where $\sqrt 5$ denotes a (fixed) square root of~5. Then
$$ \prod \left\{\,a\in \F_q^\x : \text{$2-r-a$ and $2+r+a$ are nonsquares} \,\right\} 
= 2\qquad  \text{if $q \equiv \pm1 \pmod{20}$,}  $$
$$\prod \left\{\, a \in \F_q^\x : \text{$2-r-a$ and $2+r+a$ are nonzero squares} \,\right\} 
= - \jacobi{-1}q r \qquad \text{otherwise.}$$

\bigskip\noindent{\bf Generalization of structure theorem.}  
Fix $0\ne c \in \F_q^\x$, and define a $c$-orbit to be 
$\calO_{v,c} = \left\{\, \pm v, \pm c/v \, \right\}$, where $v\in\cj\F_q^\x$. Let $f(x) = (x+c/x)^2/4$. Then for $v,w \in \cj\F_q^\x$,
$f(v) = f(w) \iff w\in\calO_{v,c} \iff \calO_{v,c} = \calO_{w,c}$. Also, if $\tau = f(v)$ then $\tau \in \F_q \iff v^q \in \calO_{v,c} \iff v^{q-1}=\pm1$ or
$v^{q+1}=\pm c$.  If $\tau  \in \F_q$, then the $c$-orbit corresponding to $\tau $ is $\left\{ \, \pm\sqrt{\tau\,} \pm \sqrt{\tau-c\,}\,\right \}$. 
If $\jacobi\tau q = A$
and $\jacobi {\tau-c}q = B$, and if $\tau \not \in \left\{\,0,c\,\right\}$, then every $v$ in the corresponding $c$-orbit 
satisfies $v^{q-AB} = A c^{(1-AB)/2}$.

The generalized structure theorem is only mildly useful, because it can be obtained from the usual one by rescaling.  Namely, if $\calO_{v,c} \longleftrightarrow \tau$
(\ie, $\tau = (v+c/v)^2/4$), then $c^{-1}\tau = (w+1/w)^2/4$, where $w=v/\sqrt c $.

\bigskip\noindent{\bf Acknowledgement.} The author thanks Dr. Art Drisko for his careful review of this article. His comments have benefitted the exposition.

\end{document}